%% file: chevalley2011dinai.tex
\documentclass[11pt,a4paper]{amsart}
\usepackage{pdfsync}

\input{oren_macros.055}

\title{Diameters of Chevalley groups\\ over local rings}
\author{Oren Dinai\\ \\ \ETHZ}
\date{\today}
\address{\aETHZ}
\email{\eETHZ /\gmail}
\begin{document}

\DO{abstract}{ Let $G$ be a Chevalley group scheme of rank $l$.  We
  show that the following holds for some absolute constant $d>0$ and
  two functions $p_0=p_0(l)$ and $C=C(l,p)$.  Let $p\ge p_0$ be a
  prime number and let $G_n:=G(\Z/p^n\Z)$ be the family of finite
  groups for $n\in \N$.

Then for any $n\ge 1$ and any subset $S$ which generates $G_n$ we have,
$$\diam(G_n,S)\le C n^d,$$
 i.e., any element of $G_n$ is a product of $C n^d$ elements from $S\cup S^{-1}$.
\IP, for some $C'=C'(l,p)$ and for any $n\ge 1$ we have, $$\diam(G_n,S)\le C' \log^d(|G_n|).$$

Our proof is elementary and effective, in the sense that the constant
$d$ and the functions $p_0(l)$ and $C(l,p)$ are calculated
explicitly. Moreover, there exists an efficient algorithm to compute a
short path between any two vertices in any Cayley graph of the groups
$G_n$.  }

\maketitle

\sec{Introduction}\l{sec:intro}

We start by recalling a few essential definitions and background
results. Let $G$ be a any group and let $S\subset G\setminus \set{1}$
be a non-empty subset.  Define $Cay(G,S)$, the (left) \emph{Cayley
  graph} of $G$ with respect to $S$, to be the \emph{undirected} graph
with vertex set $V:=G$ and edges $E:=\{ \{ g,\, sg\} :g\in G,s\in
S\}$.
\par
Now, given any finite graph $\Gamma=(V,E)$, one defines
$\diam(\Gamma)$, the diameter of $\Gamma$, to be the minimal $l\geq 0$
such that \emph{any} two vertices are connected by a path in $G$
involving at most $l$ edges (with $\diam(\Gamma)=\infty$ if the graph
is not connected). Now define $\diam(G,S)$, the diameter of a group
$G$ with respect to $S\subset G$, to be the minimal number $k$ for
which any element in $G$ can be written as a product of at most $k$
elements in $S \cup S^{-1}$.

One is naturally interested in minimizing the diameter of a group with
respect to \ti{arbitrary} set of generators. For this we define,
$$\diam(G):=\max\{\diam(G,S): S\( G \tx{ and } S
\tx{ generates } G\}.$$

\par

The diameter of groups, aside from being a fascinating field of
research, has huge amount of applications to other important
fields. In addition to Group theory and Combinatorics, the diameter of
groups in widely known for its role in Theoretical Computer Science
areas such as Communication Networks, Algorithms and Complexity (for a
detailed review about these aspects, see
\cite{babai1990diameter}). The wide spectrum of applications involved
makes this an interdisciplinary field.

It turns out that quite a lot is known about the ``best'' generators,
i.e. that a small number of well-chosen generators can produce a
relatively small diameter (see \cite{babai1990diameter}). But very
little was known until recently about the worst case. A well known
conjecture of Babai (cf. \cite{babai1988diameter,babai1992diameter})
asserts:

\begin{conj}[Babai]\label{conj babai}
  There exist two constants $d,C>0$ such that for any finite
  non-abelian simple group $G$ we have \[ \diam(G)\leq
  C\cdot\log^{d}(\left|G\right|).\] This bound may even be true for
  $d=2$, but not for smaller $d$, as the groups $Alt(n)$ demonstrate.
\end{conj}

For these type of groups, there has been enormous progress recently,
due in particular to Pyber-Szab\'o~\cite{pyber2010growth2} and
Breuillard-Green-Tao~\cite{breuillard2010linear}, when many families
of Cayley graphs of finite groups of Lie type have been shown to be
expander families (see also
\cite{helfgott2005growth,bourgain2006new,dinai2011growth} for previous
results).  However, although most of the known results are effective,
in the sense that the constants can be computed in principle, they are
usually not explicit: no specific values are given, the exception
being~\cite{kowalski2012explicit} which contains an explicit version of
Helfgott's solution of Babai's conjecture for
$\mathrm{SL}_2(\BBZ/p\BBZ)$.  But even this does not give \emph{an
  efficient algorithm} for computing a short path between any two
vertices in the Cayley graph, whose existence is guaranteed by the
diameter bounds.

In section \S\ref{sec:preliminaries} we introduce the required
definitions to be used in the next sections.  In section
\S\ref{sec:main} we are proving the main results of this manuscript,
which is \rCORO{coro_main}: this gives explicit bounds for the
constant $d$ and the functions $p_0=p_0(l)$ and $C=C(l,p)$ as stated
in the abstract.  The following is a special case of this corollary
(the precise version gives a specific value of the constant $C$).

\DO{thm}{
Let $G$ be a Chevalley group scheme of rank $l$ and dimension $k$.
Fix a prime number $p$ with $p> \max\set{\frac{l+2}2,19}$ .
Denote $G_n:=G(\Z/p^n\Z)$ for $n\in \N$.
Then any $n\ge 1$ we have,
$$\diam(G_n)\le C p^{2k} n^{10},$$ 
for some constant $C$ which depends on $G$ but not on $p$.
}

%For any $i\ge 2$ set $C_i(p,k):=\diam(G_i)$ and $d_i=d_i(3)$ where $d_i(r):=\frac{\log(4r)}{\log(2i)-\log(i+1)}$.

Although, for a fixed generating set, one can now often prove that the
relevant Cayley graphs form an expander, which provides asymptotically
a better bound, these are not usually explicit. There is also some
interest in polylogarithmic bounds for the diameter of groups:
in~\cite{ellenberg2010Expander}, there are applications of such bounds to
questions in arithmetic geometry, and there is a possibility that
explicit bounds as we have obtained could be useful to obtain more
quantitative versions of certain of those results.
\par
In section \S\ref{sec:algorithm} we explain the variant of the
``Solovay-Kitaev'' algorithm that provides fast computations of
representations of a given element as a short word, with respect to an
arbitrary set of generators.

%Let $G$ be a finite group and let $S\subset G\setminus \set{1}$ be a non-empty subset.
%Define the \emph{Cayley graph} of $G$ with respect to $S$ to be the \emph{undirected}
%graph $Cay(G,S):=\left(V,E\right)$ such that $V:=G$ and $E:=\{ \{ g,\, sg\} :g\in G,s\in S\}$.
%Given a finite graph $\Gamma=(V,E)$, define $\diam(\Gamma)$, the diameter of $\Gamma$, to be the
%minimal $l$ such that \emph{any} two vertices are connected by a path
%with at most $l$ edges. Set $\diam(\Gamma)=\infty$ if the graph is not connected.
%Now define $\diam(G,S)$, the diameter of a group $G$ with respect to $S$,
% to be $\diam(Cay(G,S))$ i.e., it is the minimal number $k$ for which any element in $G$ can be written 
% as a product of at most $k$ elements in $S \cup S^{-1}$.
%
%We will be interested in minimizing the diameter of a group with respect to \te{arbitrary} set of generators. For this we define,
%\[\diam(G):=\max\{\diam(G,S): S\subset G \text{ and } S
%\tr{ generates } G\}.\]
%The diameter of groups, aside from being a fascinating field of research, has huge amount of applications to other
%important fields. In addition to Group theory and Combinatorics, the diameter of groups in widely known for its role
%in Theoretical Computer Science areas such as Communication Networks, Algorithms and Complexity. For a detailed review see \cite{babai1990diameter}. The wide spectrum of applications of the diameter of groups, makes it an interdisciplinary and reach field.
%

\sec{Preliminaries}\l{sec:preliminaries}

First, we begin with a few preliminary definitions.

%Now We start by recalling a few essential definitions and background results.
%From this section we will allow ourself to use more technical definitions.

\doDEF{
Let $A,B$ be subsets of a group $G$ and $r\in \N$. Denote: \doitem{
  \item $A\cdot B= \set{ab:a\in A,b\in B}$.
  \item $A^{(r)}$ the subset of products of $r$ elements of $A$ with $A^{(0)}=\set{1}$.
  \item $A^{[r]}$ the subset of products of $r$ elements of $A\cup A^{-1}\cup \set{1}$.
	}
Denote the commutator word $\set{a,b}:=(ba)^{-1}ab$ and, \doitem{
	\item $\set{A,B}_1:=\set{\set{a,b}:a\in A,b\in B}$.
	\item $\set{A,B}_r$ the subset of products of $r$ elements of $\set{A,B}_1$.
	}
	
        The group $G$ will be called \ti{$r$-strongly perfect} if
        $G=\set{G,G}_r$.  Similarly if $L$ is a Lie algebra with Lie
        bracket $[a,b]$ then we replace the previous notations by
        $[A,B]_r$ and the product by summation, and $L$ will be called
        $r$-strongly perfect if $L=[L,L]_r$.  }

      \doDEF{\label{def-general} Let $G$ be a Chevalley group
        scheme\footnote{I.e., for some absolute $n\ge 1$ and for any
          commutative ring $R$ with a unit, $G(R)\le GL_n(R)$ (and
          $L(R)\le gl_n(R)$). Moreover $G$ and $L$ are functors, i.e.,
          they transform homomorphisms between objects.} associated
        with a connected complex semi-simple Lie group $G_c$ and let
        $L$ be its Lie algebra (cf. \cite{abe1969chevalley}).  Let $p$
        be a prime number and $\Z_p$ be the $p$-adic integers.  Set
        $\Gamma_0:=G(\Z_p)$, $L_0:=L(\Z_p)$ and denote for $n\ge 1$:
        \doitem{
	\item $G_n:=G(\Z_p/p^n\Z_p)\cong G(\Z/p^n\Z)$.
	\item $\pi_n$ the natural projection from $\Gamma_0$ onto
          $G_n$.
	\item $\Gamma_n:=\Gamma(p^n)=\Ker(\pi_n)$.
	\item Given $g,h\in \Gamma_0$ denote $g\EQV[n] h$ if $\pi_n(g)=\pi_n(h)$.
	\item $\Delta_n:=\Gamma_n/\Gamma_{n+1}$.
	}

        Both $\Gamma_0$ and $L_0$ have an operator ultra-metric which
        is induced by the $l_{\infty}$-norm and the absolute value on
        $\Z_p$ (which is defined, say, by $|p|=\frac12$ and then
        extended uniquely to $\Z_p$).  }

We will use the following proposition due to Weigel
\cite[Prop.~4.9]{weigel2000rigidity}. The proof for the classical
groups is easy so we give here an elementary proof of it.

\doPROP{[Weigel]\label{lem-exp-onto}
Let $G$ be a Chevalley group over $\Z_p$ and $L_0$ and $\Gamma_n$ be as in definition \ref{def-general}. Then $$\Gamma_n = \exp(p^n L_0).$$
}

\doPROOF{

  The direction $\exp(p^n L_0) \( \Gamma_n$ is trivial so we will
  prove the other direction.  We will prove only $\Gamma_1\( \exp(p
  L_0)$ since the case $n>1$ follows by the same argument.  Let $g\in
  \Gamma_1$ be $g=I+pA$ for some p-adic matrix $A$.  Since the
  summation $\ln(g)=pA-\frac12(pA)^2+\frac13(pA)^3-\ldots$ converges
  we are left to show that $\ol{\ln}(g)\in L_0$ where
  $\ol{\ln}(g):=A-\frac12 p A^2+\frac13 p^2 A^3-\ldots$ is the
  ``normalized'' logarithm.

  We can assume that $L$ is a simple Lie algebra since the statement
  holds for semi-simple Lie algebras if it holds for simple Lie
  algebras. We will prove this claim when $G$ is a classical Chevalley group i.e., of
  type $A_l,B_l,C_l$ or $D_l$.  In all these cases we will use the
  classical faithful matrix representations of $G$ and $L$ (over
  $\ol{\Q}_p$).  If $G$ is of type $A_l$ then $g\in G(\Z_p)\=
  \det(g)=1$, and $A\in L(\Z_p) \= \Tr(g)=0$.  Since
  $p\Tr(\ol{\ln}(g))=\Tr(\ln(g))=\ln(\det(g))=0$ we are
  done\footnote{We used the identity $\det(e^A)=e^{\Tr(A)}$ which is
    valid over any valuation ring (using the Jordan decomposition of
    $A$ over an algebraic close field extending the ring).} in this
  case.

  Now suppose $G$ is a Chevalley group of type $B_l,C_l$ or
  $D_l$. Then we have a vector space $V$ of finite dimension (over
  $\Q_p$) with some non-singular bi-linear form $\gb$ on $V$. For
  $A\in \End(V)$ denote by $A^*$ the $\gb$-adjoint\footnote{So that
    $A\mapsto A^*$ is an anti-automorphism of $\End(V)$ of order $2$
    with $\gb(Av,w)\equiv \gb(v,A^*w)$.} of $A$.  Then $g\in G(\Z_p)\=
  gg^*=I$, and $A\in L(\Z_p) \= A+A^*=0$.  Since $\ln(g)$ and
  $\ln(g^*)=\ln(g)^*$ converge and $g,g^*$ commute we get that
 $$\ln(gg^*)=\ln(g)+\ln(g)^*=p(\ol{\ln}(g)+(\ol{\ln}(g))^*)=\ln(I)=0,$$ so we are done in these cases as well.

}

\doDEF{\label{def_diam_chain}
Let $N\le H \le G$ be a chain of groups (not necessarily normal) and $S\( G$. Denote:
	\doitem{
		\item $\diam(H/N;S)=\min\set{l:H \( S^{[l]}N}$.
		\item $\diam_G(H/N):=\max\set{\diam(H/N;S):\gen{S}=G}$.
		\item $\diam(H/N):=	\diam_H(H/N)$.
	}
}

Note that $\diam(H/N)$ is the worst diameter of the Schreier graphs of
$H/N$ and if $N=1$ then this is the worst diameter of the Cayley
graphs of $H$.

\doSFACT{\label{fact_diam_chain}
Let $N\le H \le G$ be a chain of groups and $S\( G$. Then,
\doitem{
	\item $\diam(G/N;S)\le \diam(G/H;S)+\diam(H/N;S)$.
	\item $\diam(G/N)\le \diam_G(G/H)+\diam_G(H/N)$.
	}
}

\sec{Main results}\l{sec:main}

\lTHM{thm_perfect_and_commutators}{
Suppose $L(\Z_p)$ is $r$-strongly perfect. Then for any $i,j\in\N$, $$\Delta_{i+j}=\set{\Delta_i,\Delta_j}_r.$$
}

\doPROOF{
The direction $[\)]$:
This is clear since $\set{\Gamma_i,\Gamma_j}_r \( \Gamma_{i+j}$.
Moreover, if $g,g'\in \Gamma_0$ and $g\EQV[i+1]I+p^i A$, $g'\EQV[j+1]I+p^j A'$ for some matrices $A,A'$, then $\set{g,g'}\EQV[i+j+1] I+p^{i+j} [A,A']$.

The direction [$\($]:
Let $g\in \Gamma_n/\Gamma_{n+1}$ with $n=i+j$.
By \rLEM{lem-exp-onto}, $g\us[n+1]{\equiv} \exp(p^n A)$ for some $A\in L_0$. \TF\ $g\EQV[n+1] I+p^n A$.
By the assumption, $A=\sum_{k=1}^r [A_k,A_k']$ for some $A_1,A_1',\ldots,A_r,A_r'\in L_0$.
Denote $g_k:=\exp(p^i A_k)\in \Gamma_i$ and  $g'_k:=\exp(p^j A'_k)\in \Gamma_j$.
Therefore $g_k \us[i+1]{\equiv}I+p^i A_k$ and $g'_k \us[j+1]{\equiv}I+p^j A_k$ and
$$g\US[n+1]{\equiv} I+p^n A \us[n+1]{\equiv} \set{g_1,g'_1}\cdot\ldots \set{g_r,g'_r}.$$
}

%If $l=2m$ is even then take $h=1\del_1+2\del_3+\ldots m\del_{l-1}$ and if $l=2m+1$ is odd then take
%$h=1\del_1+2\del_3+\ldots m\del_l$
%Define $h=1\eps_1-1\eps2+2\eps_3-2\eps_4\ldots $

\lLEM{lem_chevalley_strongly_perfect}{
Let $G$ be a Chevalley group of rank $l$, $L$ its Lie algebra and let $p\ge \frac{l+2}2$ be an odd prime number.
If $G$ is a group of exceptional Lie type then suppose that $p>19$.
Then $L(\Z_p)$ is $3$-strongly perfect.
}

\doPROOF{

  Let $B=\set{e_s,h_r: s\in \Phi, r \in \Pi}$ be a Chevalley basis of
  $L$, where $\Phi$ is the root system associated to $L$ and $\Pi$ are
  the simple roots of $\Phi^+$ (for some fixed order).
  \WLOG,\footnote{Since the statement holds for semi-simple Lie
    algebras if it holds for simple Lie algebras.} we can assume that
  $\Phi$ is irreducible.

%Then $l=\rank(L)=|\Pi|$ and $k=\dim(L)=|\Pi|+|\Phi|$.
For any $r\in \Phi$ denote $L_r:=\Z_p e_r$ and $H_r:=\Z_p h_r$ where $h_r=[e_r,e_{-r}]$ is the co-root of $r$.
We have $L(\Z_p)=L_{\Phi}\oplus H$ where $H:=\bigoplus_{r\in \Pi} H_r$ and $L_{\Phi}:=\bigoplus_{r\in \Phi} L_r$.
We will use the following facts about Lie bracket of the root system.
For any $h\in H$ and $s\in \Phi$ we have $[h,e_s]=(h,s)e_s$ where $(\cdot,\cdot)$ is the inner product in $H$.
For any linearly independent pair of roots (i.e., $r\neq \pm s$) we have $[e_r,e_s]\in L_{\Phi}$ and if their sum $r+s\notin \Phi$ then $[e_r,e_s]=0$.

We will say that a submodule $V\le L(\Z_p)$ is \ti{covered} if $V\( [L(\Z_p),L(\Z_p)]$.
For any $X\( \Phi$ denote $L_X:=\bigoplus_{r\in X} L_r$. We will say that $X$ is \ti{covered} if there exists $h\in H$ with $(h,X)\( (\Z_p)^{\times}$.
We will say that $\Phi$ is \ti{$k$-covered} if $\Phi= X_1\cup  \ldots \cup  X_k$ and each $X_i$ is covered.
Note that if $X$ is covered by some $h$ then $L_X$ is covered; indeed if $y=\sum a_r e_r\in L_X$ then $[h,y']=y$
where $y'= \sum \frac{a_r}{(r,h)} e_r \in L_X$.

Note also that $H$ is always covered; indeed for any $x=\sum a_r h_r \in H$ we have $x=[x',x'']$ where
$x'=\sum_{r\in \Pi} a_r e_r$ and $x''=\sum_{r\in \Pi} e_{-r}$.
In order to complete the proof we will show that $\Phi$ is $2$-covered.

We will use the following notations. Suppose that $\Phi$ can be embedded into an Euclidean space $E\cong H$ of
dimension $l$ \ST\ $\set{\ga_i}$ is an orthonormal basis of $E$.
Set $h_1:=\sum \ga_i \in H$ and $h_2:=\sum \lam_i \ga_i \in H$ where $\lam_1,\ldots,\lam_l \in \Z \cap (-p,p)$
and for any $i\neq j$ we have $\lam_i-\lam_j\in \Z \SM p\Z$; e.g., we can take the $\lam_i$'s to be a subset
of $\set{0,\pm 1, \pm 2,\ldots,\pm \frac{p-1}2}$. Later we will put more restrictions on the choice of the $\lam_i$'s.

First suppose that $\Phi$ is one of the classical root systems.
If $\Phi = A_l$ then by \cite{dinai2006poly} it is $2$-strongly perfect since $\Phi$ is covered (cf. \cite{gamburd2004uniform}).
Now suppose $\Phi$ is of type $B_l,C_l$ or $D_l$.
Set $\Phi= X_1\cup  X_2$ where $X_1\( \set{\pm (\ga_i-\ga_j):i\neq j}$ and $X_2\( \set{\pm (\ga_i+\ga_j),\pm \ga_i, \pm 2\ga_i:i\neq j}$.
If $p>2$ then $(h_1,X_1)\( \set{\pm1,\pm2}\( (\Z_p)^{\times}$.
If in addition $2(p-1)\ge l$ then we can find $\lam_1,\ldots,\lam_l$ as above \ST\ $\sum \lam_i =0$; \tF\ $(h_2,X_2)\( (\Z_p)^{\times}$.
We got that the classical root systems are $2$-covered and so they are $3$-strongly perfect.

Now we shall see that essentially the same argument works if $\Phi$ is an exceptional root systems (cf. \cite[\S8]{carter2005lie} for a complete list of roots of each type). If $\Phi$ is of type $G_2$ with $l=2$ then$(h_1,X_1)\( \set{\pm1,\ldots,\pm 5}$; \tF\ $(h_1,X_1)\( (\Z_p)^{\times}$ provided $p>5$; so $\Phi$ is $1$-covered provided $p\ge 5$.

Now suppose $\Phi$ is of type $F_4$ and $\Phi=\set{\pm \ga_i,\pm \ga_i\pm \ga_j,\sum_{k=1}^4\pm \ga_k:i\neq j}$.
Split the set $\Phi= X_1\cup  X_2$ where $X_1$ is the ``unbalanced'' subset of sums where the number of $+$'s is not equal to
the number of $-$'s and $X_2$ is the ``balanced'' subset; i.e., $X_2:=\set{\ga_i-\ga_j,\ga_{i_1}+\ga_{i_2}-\ga_{i_3}-\ga_{i_4}}$.
Set $\set{\lam_i}=\set{0,1,\pm2}$. Then $(h_1,X_1)\( \set{\pm1,\pm2,\pm4}$ and $(h_2,X_2)\( \set{\pm1,\pm2,\pm4}$.
\TF\ $\Phi$ is $2$-covered provided $p\ge 5$.

Now let us show that $E_8$ is $3$-covered (and therefore also $E_6,E_7$).
Now $l=8$ and again we split $\Phi$ into an unbalanced set $X_1$ and a balanced set $X_2$.
Set $\set{\lam_i}=\set{0,1,\pm2,\pm3,\pm4}$. Then $(h_1,X_1)\( \pm \set{2,4,8}$ and $(h_2,X_2)\( \pm \set{1,\ldots,19}$.
\TF\ we get that $\Phi$ is $2$-covered provided $p>19$; so we are done.
}

\lCORO{coro_main}{
Let $G$ be a Chevalley group of rank $l$ and let $p$ be a prime number chosen as above.
Denote $G_n:=G(\Z/p^n\Z)$ for $n\in \N$.
For any $i\ge 2$ set $C_i(p,k):=\diam(G_i)$ and $d_i=d_i(3)$ where $d_i(r):=\frac{\log(4r)}{\log(2i)-\log(i+1)}$.
Then for any $n\ge 1$ and $i\ge 2$ we have\footnote{\IP\ $C_i\le p^{ik}$ where $k=\dim(L)$ and $d_i$ is monotone decreasing to $2+\log_2(3)$.},
$$\diam(G_n)\le C_i n^{1+d_i}.$$
}

\doPROOF{
Denote $L_n(j)=\diam_{G_n}(\Delta_j)$ for $0\le j< n$. Then by \rFACT{fact_diam_chain},
$$\diam(G_n)\le L_n(0)+L_n(1)+\ldots+L_n(n-1).$$
%$$\diam(G_n)\le \diam_{G_n}(\Delta_0)+\diam_{G_n}(\Delta_1)+\ldots \diam_{G_n}(\Delta_{n-1}).$$
By induction on $j$, we will prove that for any $i\ge 2$ and $0\le j< n$, $$L_n(j)\le C_i j^{d_i},$$ and therefore,
$$\diam(G_n)\leq \sum_{j=0}^{n-1} C_i j^{d_i}\le C_i n^{1+d_i},$$
as we claimed.

Fix some $i\ge 2$. The induction base is for $j<i$, and then trivially $L_n(j)\le \diam(G_i)=C_i$.
Now suppose $j\ge i$. Then by \rTHM{thm_perfect_and_commutators}, by \rLEM{lem_chevalley_strongly_perfect} with $r=4$ and by the induction assumption, we get
$$L_n(j)\le 4r L_n(\lfloor \frac{j+1}2\rfloor)\le 4r C_i (\frac{j+1}2)^{d_i}=4r (\frac{j+1}{2j})^{d_i}C_i j^{d_i}\le C_i j^{d_i},$$
since by the definition of $d_i$, $4r (\frac{j+1}{2j})^{d_i}\le 1$ for any $j\ge i$.

}

\DO{remark}{ The combination of \rTHM{thm_perfect_and_commutators},
  \rLEM{lem_chevalley_strongly_perfect} and \rCORO{coro_main} give a
  generalization of what is known as the ``Solovay-Kitaev method''.

  Geometrically we divide the group $\Gamma_0$ into neighborhoods of
  the identity $\Gamma_n$, and their ``layers'' $\Delta_n$.  First, we
  use the global properties of the Lie brackets in order to get local
  properties of the commutators in these layers.  Then
  \rCORO{coro_main} allows us to ``glue'' the local properties valid
  in these layers into a global property.

  Note that this method can prove, at best, a bound of order of
  magnitude $\log^d(|G|)$, with $d$ arbitrary close to $2$, but not a
  better bound. This follows because the best possible situation is
  that $L$ is $1$-strongly perfect.  }

\sec{The Solovay-Kitaev algorithm}\l{sec:algorithm}

Now we give an explicit description and analysis of the Solovay-Kitaev
algorithm (cf. \cite[\S 3]{dawson2005solovay} and also
\cite{nielsen2002quantum}).  First we describe a procedure based on
\rTHM{thm_perfect_and_commutators} and
\rLEM{lem_chevalley_strongly_perfect} from the previous section. This
procedure is an effective version of these statements about finding an
explicit decomposition of an element as a product of (at most four)
commutators.

\ssec{Commutator decomposition}

The main algorithm (in the next section) will use the subalgorithm
$SK'(g,n),$ which gets an input $g\in \Gamma_n$ with $n\ge 2$; then it
returns a pair of quadruples $((g_i),(g_i'))$ such that
$\set{g_1,g_1'}\cdot\ldots \set{g_4,g_4'}\US[n+1]{\equiv} g$ where
$g_i,g_i'\in \Gamma_m$ with $m\ge \frac{n-1}2$. Note that this is a
direct consequence \rTHM{thm_perfect_and_commutators} and
\rLEM{lem_chevalley_strongly_perfect}; if $g\us[n+1]{\equiv} \exp(p^n
A)\EQV[n+1] I+p^n A$ for some $A\in L_0$ and $A=\sum_{k=1}^r
[A_k,A_k']$ (with $r=4$) then by \rTHM{thm_perfect_and_commutators} we
get the required solution $g\US[n+1]{\equiv} \set{g_1,g'_1}\cdot\ldots
\set{g_r,g'_r}$; in order to solve $A=\sum_{k=1}^r [A_k,A_k']$ we
first find the decomposition of $A$ as a linear combination in the
Chevalley basis and then use \rLEM{lem_chevalley_strongly_perfect} in
order to decompose it as a sum of (at most) four Lie brackets.

\ssec{The Solovay-Kitaev algorithm}

The Solovay-Kitaev algorithm $SK(g,\ol{s},n)$ gets an element $g\in \Gamma_0$, $n\in\N$ and a $m$-tuple $\ol{s}$ (with entries in $\Gamma_0$) that generates $G_n=\Gamma_0/\Gamma_n$; then it returns a word $w\in F_m$ (in $m$ letters) such that $g\us[n]{\equiv} w(\ol{s})$.
If $n \le 2$ then $SK$ returns such a word simply by checking all the possible words of length $l(w) \le |G_2|=|G(\Z/p^2\Z)|$.
If $n>2$, set $w_0=SK(g,\ol{s},n-1)$ and $z=w_0(\ol{s})^{-1}g\in \Gamma_{n-1}$ and let $(\ol{x},\ol{y})=SK'(z,n-1)$.
Set for $k=1,\ldots 4$, $w_k:=SK(x_k,\ol{s},n-1)$ and $w_k':=SK(y_k,\ol{s},n-1)$ and return $w:=w_0\cdot \set{w_1,w_1'}\cdot\ldots \set{w_4,w_4'}$.

\ssec{Analysis of the algorithm}

The return length of the output word of the algorithm is $C_i n^{1+d_i}$, the same as was described in \rCORO{coro_main}.
Note that $d_2<9$; $C_i\le p^{ik}$ where $k=\dim(L)=|\Phi|+|\Pi|$ ; and $d_i$ is monotone decreasing to $2+\log_2(3)$.

\begin{acknowledgement*}
I would like to thank Alex Lubotzky for bringing Thomas Weigel's results to my attention and Emmanuel Kowalski for 
many helpful comments and suggestions.
\end{acknowledgement*}

%\XX{add time/space\\ complexity?}

%\cite{dinai2006uniform,dinai2009growth,dinai2010expansion}

\useBIB{amsalpha}{my_bib.055}
%\useBIB{bib_style_wdg_jmc}{my_bib.003}

\end{document}

%% file: oren_macros.055.tex
\usepackage{amsmath}
\usepackage{amsthm}
\usepackage{amssymb}
\usepackage{amsfonts}

\def\ncomm{\newcommand}
\def\rncomm{\renewcommand}
\def\NCOM{\newcommand}

\def\nthm{\newtheorem}
\def\ni{\noindent}

\ncomm{\doit}[2]{\begin{#1} #2 \end{#1}}
\def\DO{\doit}

\ncomm{\MAX}[2][{}]{\max_{#1}\!\set{#2}}
\ncomm{\MIN}[2][{}]{\min_{#1}\!\set{#2}}
\ncomm{\CUP}[1][{}]{\bigcup_{#1}}
\ncomm{\CAP}[1][{}]{\bigcap_{#1}}

\ncomm{\COP}[1][{}]{\coprod_{#1}}

\ncomm{\mbold}[1]{\boldsymbol{#1}}
\ncomm{\mbf}[1]{\em{\boldsymbol{#1}}}

\theoremstyle{plain}
\nthm{thm}{Theorem}[section]
\nthm{lemma_}[thm]{Lemma}
\nthm{fact_}[thm]{Fact}
\nthm{prop}[thm]{Proposition}
\nthm{corr}[thm]{Corollary}
\nthm{coro}[thm]{Corollary}
\nthm{assumption}[thm]{Assumption}
\nthm{example}[thm]{Example}
\nthm{examples}[thm]{Examples}
\nthm{remark}[thm]{Remark}
\nthm{remarks}[thm]{Remarks}
\nthm{question}[thm]{Question}

\theoremstyle{definition}
\nthm*{notation}{Notation}
\nthm{quest}{Question} 
\nthm{conj}{Conjecture}
\nthm*{algorithm}{Algorithm}

\nthm{def_}[thm]{Definition}
\nthm{sfact}[thm]{Simple Fact}

\theoremstyle{remark}
\nthm*{rem}{Remark}
\nthm*{mynote}{Note}
\nthm{case}{Case}

\nthm*{acknowledgement*}{Acknowledgement} 
\nthm*{thanks*}{Acknowledgement}

\ncomm{\addrem}[1]{\doit{rem}{#1}}
\ncomm{\addnote}[1]{\doit{mynote}{#1}}
\ncomm{\addcase}[1]{\doit{case}{#1}}

\ncomm{\doproof}[1]{\doit{proof}{#1}}

\ncomm{\dothm}[1]{\doit{thm}{#1}}
\ncomm{\docorr}[1]{\doit{corr}{#1}}
\ncomm{\docoro}[1]{\doit{corr}{#1}}
\ncomm{\dofact}[1]{\doit{fact}{#1}}
\ncomm{\dosfact}[1]{\doit{sfact}{#1}}
\ncomm{\dodef}[1]{\doit{def_}{#1}}
\ncomm{\doTHANKS}[1]{\doit{thanks*}{#1}}

\ncomm{\doquest}[1]{\doit{quest}{#1}}
\ncomm{\doprob}[1]{\doit{pro}{#1}}
\ncomm{\doconj}[1]{\doit{conj}{#1}}

\ncomm{\doprop}[1]{\doit{prop}{#1}}
\ncomm{\dolemm}[1]{\doit{lemma_}{#1}}
\ncomm{\dolem}[1]{\doit{lemma_}{#1}}
\ncomm{\donota}[1]{\doit{notation}{#1}}

\ncomm{\doTHM}[1]{\dothm{#1}}
\ncomm{\doCORO}[1]{\docoro{#1}}
\ncomm{\doPROP}[1]{\doprop{#1}}
\ncomm{\doLEM}[1]{\dolem{#1}}
\ncomm{\doFACT}[1]{\dofact{#1}}
\ncomm{\doSFACT}[1]{\dosfact{#1}}
\ncomm{\doDEF}[1]{\dodef{#1}}

\ncomm{\doEX}[1]{\doit{example}{#1}}
\ncomm{\doEXS}[1]{\doit{examples}{#1}}
\ncomm{\doREM}[1]{\doit{remark}{#1}}
\ncomm{\doREMS}[1]{\doit{remarks}{#1}}

\ncomm{\doPROOF}[1]{\doproof{#1}}

\ncomm{\lTHM}[2]{\dothm{\l{#1}#2}}
\ncomm{\lCORO}[2]{\docoro{\l{#1}#2}}
\ncomm{\lPROP}[2]{\doprop{\l{#1}#2}}
\ncomm{\lLEM}[2]{\dolem{\l{#1}#2}}
\ncomm{\lFACT}[2]{\dofact{\l{#1}#2}}
\ncomm{\lSFACT}[2]{\dosfact{\l{#1}#2}}
\ncomm{\lDEF}[2]{\dodef{\l{#1}#2}}

\ncomm{\prTHM}[2]{\doproof{[Proof of Theorem \r{#1}]#2}}
\ncomm{\prCORO}[2]{\doproof{[Proof of Corollary  \r{#1}]#2}}
\ncomm{\prPROP}[2]{\doproof{[Proof of Proposition \r{#1}]#2}}
\ncomm{\prLEM}[2]{\doproof{[Proof of Lemma \r{#1}]#2}}
\ncomm{\prFACT}[2]{\doproof{[Proof of Fact \r{#1}]#2}}
\ncomm{\prSFACT}[2]{\doproof{[Proof of Simple Fact \r{#1}]#2}}

\ncomm{\rTHM}[1]{Theorem \r{#1}}
\ncomm{\rCORO}[1]{Corollary  \r{#1}}
\ncomm{\rPROP}[1]{Proposition \r{#1}}
\ncomm{\rLEM}[1]{Lemma \r{#1}}
\ncomm{\rFACT}[1]{Fact \r{#1}}
\ncomm{\rSFACT}[1]{Simple Fact \r{#1}}
\ncomm{\rDEF}[1]{Definition \r{#1}}

\ncomm{\boxedeqn}[1]{%
  \[\fbox{%
      \addtolength{\linewidth}{-2\fboxsep}%
      \addtolength{\linewidth}{-2\fboxrule}%
      \begin{minipage}{\linewidth}%
      \begin{equation}#1\end{equation}%
      \end{minipage}%
    }\]%
}

\newcommand{\eval}[2][\right]{\relax
  \ifx#1\right\relax \left.\fi#2#1\rvert}

\ncomm{\doeq}[1]{\doit{equation*}{#1}}       %%for short equations
\ncomm{\doEQ}[2]{\doit{equation}{\l{#1}#2}}

\ncomm{\mnote}[1]{\marginpar{#1}}
\ncomm{\MARG}[1]{\marginpar{#1}}
\ncomm{\XXX}[1]{{\marginpar{#1??}}XXX}
\ncomm{\XX}[1]{{\marginpar{- #1}}$^{(?)}$}
\ncomm{\xxx}[1]{{\marginpar{#1\\-\tb{xxx?}}}$^{\tb{xxx?}}$}

\ncomm{\doeqsA}[2]{\subeq{#1 \doAlign{#2}}}

\ncomm{\subeq}[1]{#1}

\ncomm{\doeqa}[2]{\doeq{#1 \doAlign{#2}}}
\ncomm{\doeqA}[2]{\doeq{#1 \doAlign{#2}}}

\ncomm{\geqBY}[1]{\byeq{#1}{\geq}}
\ncomm{\leqBY}[1]{\byeq{#1}{\leq}}
\ncomm{\gBY}[1]{\byeq{#1}{>}}
\ncomm{\lBY}[1]{\byeq{#1}{<}}
\ncomm{\neqBY}[1]{\byeq{#1}{\neq}}

\ncomm{\geBY}[1]{\byeq{#1}{\ge}}
\ncomm{\leBY}[1]{\byeq{#1}{\le}}
\ncomm{\gnBY}[1]{\byeq{#1}{>}}
\ncomm{\lnBY}[1]{\byeq{#1}{<}}
\ncomm{\neBY}[1]{\byeq{#1}{\ne}}

\ncomm{\eqBY}[1]{\byeq{#1}{=}}
\ncomm{\ggBY}[1]{\byeq{#1}{\gg}}
\ncomm{\llBY}[1]{\byeq{#1}{\ll}}
\ncomm{\subseteqBY}[1]{\byeq{#1}{\subseteq}}
\ncomm{\nsubseteqBY}[1]{\byeq{#1}{\nsubseteq}}
\ncomm{\notsubseteqBY}[1]{\byeq{#1}{\notsubseteq}}
\ncomm{\RightarrowBY}[1]{\byeq{#1}{\Longrightarrow}}
\ncomm{\LeftarrowBY}[1]{\byeq{#1}{\Longleftarrow}}

\def\>{\Rightarrow}

\def\<{\Leftarrow}

\def\={\Leftrightarrow}
\def\({\subseteq}
\def\){\supseteq}

\ncomm{\eqcenter}[1]{&\;{#1}\;&}
\ncomm{\geqBYc}[1]{\eqcenter{\geqBY{#1}}}
\ncomm{\leqBYc}[1]{\eqcenter{\leqBY{#1}}}
\ncomm{\gBYc}[1]{\eqcenter{\gBY{#1}}}
\ncomm{\lBYc}[1]{\eqcenter{\lBY{#1}}}
\ncomm{\eqBYc}[1]{\eqcenter{\eqBY{#1}}}
\ncomm{\neqBYc}[1]{\eqcenter{\neqBY{#1}}}
\ncomm{\ggBYc}[1]{\eqcenter{\ggBY{#1}}}
\ncomm{\llBYc}[1]{\eqcenter{\llBY{#1}}}
\ncomm{\subseteqBYc}[1]{\eqcenter{\subseteqBY{#1}}}
\ncomm{\nsubseteqBYc}[1]{\eqcenter{\nsubseteqBY{#1}}}
\ncomm{\notsubseteqBYc}[1]{\eqcenter{\notsubseteqBY{#1}}}
\ncomm{\RightarrowBYc}[1]{\eqcenter{\RightarrowBY{#1}}}
\ncomm{\LeftarrowBYc}[1]{\eqcenter{\LeftarrowBY{#1}}}

\ncomm{\rightarrowUP}[1]{\UP{#1}{\rightarrow}}
\ncomm{\mapstoUP}[1]{\UP{#1}{\mapsto}}
\ncomm{\dostam}[1]{}

\ncomm{\set}[1]{\left\{#1\right\}}
\ncomm{\sset}[1]{\left|\left\{#1\right\}\right|}
\ncomm{\SER}[3][1]{{#3}_{#1},\ldots,{#3}_{#2}}
\ncomm{\SET}[3][1]{\set{\SER[#1]{#2}{#3}}}
\ncomm{\VEC}[3][1]{\lr{\SER[#1]{#2}{#3}}} 

\ncomm{\dotop}[2][]{#1\atop #2}
\ncomm{\stack}[1]{\substack{#1}}
\ncomm{\stackl}[1]{\begin{subarray}{l}#1\end{subarray}}
      
\ncomm{\UP}[2][]{\stackrel{#1}{#2}}
\ncomm{\os}[2][]{\overset{#1}{#2}}
\ncomm{\us}[2][]{\underset{#1}{#2}}
\ncomm{\OS}[2][]{\overset{_{#1}}{#2}}
\ncomm{\US}[2][]{\underset{^{#1}}{#2}}
\ncomm{\EQV}[1][]{\US[#1]{\equiv}}

\ncomm{\byeq}[2]{\stackrel{\re{#1}}{#2}}
\ncomm{\BY}[2]{\stackrel{\re{#1}}{#2}}
\def\em{\ensuremath}

\ncomm{\xydiagram}[8]{
\xymatrix{ 
    {#1} \ar[r]^{#2} \ar[d]_{#4} & {#3} \ar[d]^{#5} \\
    {#6} \ar[r]_{#7}       & {#8}
}}

\ncomm{\tn}[1]{\textnormal{{#1}}}
\ncomm{\tb}[1]{\textbf{{#1}}}
\ncomm{\tu}[1]{\underline{{#1}}} 
\ncomm{\te}[1]{\emph{{#1}}}
\ncomm{\ti}[1]{\textit{{#1}}}
\ncomm{\tr}[1]{\textrm{{#1}}}
\ncomm{\tx}[1]{\text{{#1}}}
\rncomm{\tt}[1]{\texttt{{#1}}}
\ncomm{\tin}[1]{{\tiny #1}}
\ncomm{\tl}[1]{{\large #1}}
\ncomm{\tll}[1]{{\Large #1}}
\ncomm{\tL}[1]{{\LARGE #1}}

\ncomm{\tup}[1]{\textup{#1}}

\ncomm{\mT}[1]{{\textstyle #1}}
\ncomm{\mD}[1]{{\displaystyle #1}}
\ncomm{\mS}[1]{{\scriptstyle #1}}
\ncomm{\mSS}[1]{{\scriptscriptstyle #1}}
\ncomm{\mBF}[1]{\mbox{\boldmath$#1$}}

\DeclareMathOperator{\diam}{diam}
\DeclareMathOperator{\End}{End}
\DeclareMathOperator{\Ker}{Ker}

\DeclareMathOperator{\Tr}{Tr}

\def\Char{{\operatorname{char}}}
\ncomm{\dochar}[1]{\Char(#1)}
\ncomm{\doChar}[1]{\op{Char}(#1)}
\def\dim{{\operatorname{dim}}}

\ncomm{\doIm}[1]{\op{Im}(#1)}

\def\Span{{\hbox{\rm span}}}
\ncomm{\dospan}[1]{\Span(#1)}
\ncomm{\doSpan}[1]{\op{Span}(#1)}

\def\ga{\alpha}
\def\gb{\beta}

\def\BBZ {{\mathbb Z}}

\ncomm{\doOP}[1]{\defOP{#1}}
\ncomm{\doOPx}[1]{\defOPx{#1}}

\ncomm{\myOP}[1]{\operatorname{#1}}
\ncomm{\myOPx}[1]{\operatorname*{#1}}
\ncomm{\op}[1]{\operatorname{#1}}

\ncomm{\vs}[1][5]{\vskip #1mm}
\ncomm{\doenum}[1]{\begin{enumerate} #1 \end{enumerate}}
\ncomm{\doitem}[1]{\doit{itemize}{#1}}

\ncomm{\dodesc}[1]{\doit{description}{#1}}

\ncomm{\doeqn}[1]{\doit{eqnarray}{#1}}

\ncomm{\arr}[1]{\doit{array}{#1}}
\ncomm{\tavla}[2]{\begin{array}{#1} #2 \end{array}}

\ncomm{\nnum}{\nonumber}
\ncomm{\nt}{\notag}
\ncomm{\nn}{\nonumber}

\ncomm{\pageB}{\pagebreak}
\ncomm{\lineB}{\linebreak}

\ncomm{\TBD}[1]{}
\ncomm{\old}[1]{}

\ncomm{\donone}[1]{#1}
\ncomm{\lr}[1]{\left(#1\right)}
\ncomm{\lrv}[1]{\left|#1\right|}
\ncomm{\lrn}[1]{\left\|#1\right\|}
\ncomm{\lrb}[1]{\left[#1\right]}
\ncomm{\lrs}[1]{\left\{\!{#1}\!\right\}}
\ncomm{\lrp}[1]{\left(\!{#1}\!\right)}
\ncomm{\lrg}[1]{\left<\!{#1}\!\right>}
\ncomm{\lre}[1]{\left.{#1}\right.}
\ncomm{\lrf}[1]{\lfloor{#1}\rfloor}

\ncomm{\domatrix}[1]{\doit{matrix}{#1}}
\ncomm{\dosmatrix}[1]{\lr{\doit{smallmatrix}{#1}}}

\ncomm{\dobmatrix}[1]{\dopmatrix{#1}}
\ncomm{\dopmatrix}[1]{\doit{pmatrix}{#1}}
\ncomm{\dobdiag}[1]{\dopmatrix{{#1}&0\\0& #1^{-1}}}
\ncomm{\dosbdiag}[1]{\dosmatrix{{#1}&0\\0& #1^{-1}}}

\ncomm{\docases}[1]{\doit{cases}{#1}}
\ncomm{\dodarrow}[2]{\left. #1 \right \downarrow_{#2}}
\ncomm{\domid}[2]{\left. #1 \right|_{#2}}

\ncomm{\dosideset}[3][]{\sideset{#1}{#2}#3}
\ncomm{\dolim}[2]{\limits_{#1}^{#2}}

\ncomm{\xright}[2][]{\xrightarrow[#1]{#2}}
\ncomm{\xleft}[2][]{\xleftarrow[#1]{#2}}

\ncomm{\tfr}[2][1]{\tfrac{#1}{#2}}
\ncomm{\fr}[2][1]{\frac{#1}{#2}}
\ncomm{\dfr}[2][1]{#1\!/#2}
\ncomm{\sqr}[1]{\sqrt{\!\!#1}}

\def\ST{such that}

\def\WLOG{Without loss of generality}

\ncomm{\stam}[1]{}
\ncomm{\none}[1]{#1}
\ncomm{\gen}[1]{\langle #1\rangle}

\def\l{\label}

\ncomm{\tagl}[2]{\tag{#1}\label{#2}}
\ncomm{\tagrl}[3]{\tagl{\ref{#1}$#2$}{#3}}

\def\r{~\ref}
\ncomm{\er}[1]{~\eqref{#1}}
\ncomm{\re}[1]{~\eqref{#1}}

\ncomm{\cf}[1]{(cf. \cite{#1})}
\ncomm{\cfr}[2][XXX]{(cf. \cite[#1]{#2})}
\ncomm{\ci}[1]{~\cite{#1}}
\ncomm{\CITE}[3][]{\tn{(#1\cite[#2]{#3})}} %XX TBD slant

\ncomm{\doenv}[2]{\noindent \tb{{#1}:}{#2}}

\ncomm{\doearray}[1]{\doit{eqnarray*}{#1}}
\ncomm{\doeArray}[1]{\doit{eqnarray}{#1}}
\ncomm{\dm}[1]{\doit{displaymath}{#1}}

\ncomm{\marray}[1]{\doit{math}{\doit{array}{#1}}}

\ncomm{\dmarray}[1]{\doit{displaymath}{\doit{array}{#1}}}

\ncomm{\doalign}[1]{\doit{align*}{#1}}
\ncomm{\doAlign}[1]{\doit{align}{#1}}

\ncomm{\doAlignSplit}[2]{\doAlign{\l{#1}\doSplit{#2}}}
\ncomm{\doalignL}[1]{\doit{flalign*}{#1}}
\ncomm{\doAlignL}[1]{\doit{flalign}{#1}}

\ncomm{\doAlignSplitL}[2]{\doAlignL{\l{#1}\doSplit{#2}}}

\ncomm{\eqL}[1]{\lefteqn{#1}}

\def\qu{\quad}

\ncomm{\tqq}[1][]{\tx{#1}\qquad}
\ncomm{\tq}[1][]{\tx{#1}\quad}
\ncomm{\qqt}[1][]{\qquad\tx{#1}}
\ncomm{\qtq}[1][]{\quad\tx{#1}\quad}

\ncomm{\qmq}[1]{\;\qu{#1}\qu\;}

\ncomm{\doalignAT}[2][2]{\doit{alignat*}{{#1}#2}}
\ncomm{\doAlignAT}[2][2]{\doit{alignat}{{#1}#2}}

\ncomm{\tBY}[1][]{&& \tx{#1}}
\ncomm{\tBYeq}[1]{\tBY[by \re{#1}]}

\ncomm{\doSplit}[1]{\doit{split}{#1}}
\ncomm{\eqSplit}[2]{\doeq{\l{#1}\doSplit{#2}}}
\ncomm{\eqsplit}[2]{\doeqnn{\doSplit{#2}}}
\ncomm{\doGather}[1]{\doit{gather}{#1}}

\def\sec{\section}
\def\ssec{\subsection}

\def\IP{In particular}
\def\TF{Therefore}
\def\tF{therefore}

\ncomm{\epss}[1][]{\varepsilon\!^{#1}}
\ncomm{\Oeps}[1][]{O(\epss[#1]\,\!)}

\ncomm{\NMID}[1]{\!\!\nmid_{#1}}
\ncomm{\MID}[1]{|_{#1}}

\def\WLOG{Without loss of generality}

\def\Q{\mathbb{Q}}

\def\N{\mathbb{N}}

\def\Z{\mathbb{Z}}

\def\gmail{oren.dinai@gmail.com}

\def\eETHZ{oren.dinai@math.ethz.ch}
\def\ETHZ{ETH Zurich}
\def\aETHZ{Oren Dinai, Department of Mathematics - ETH Zurich, Ramistrasse 101, 8092 Zurich, Switzerland.}

\ncomm{\seq}[3][1]{#3_{#1},\ldots,\,#3_{#2}}

\ncomm{\n}[1]{\em{\|#1\|}}
\ncomm{\Lo}[1]{\em{\n{#1}_{1}}}
\ncomm{\Lt}[1]{\em{{\n{#1}}_2}}
\ncomm{\Lp}[1]{\em{\n{#1}_p}}
\ncomm{\Lq}[1]{\em{\n{#1}_q}}
\ncomm{\Lr}[1]{\em{\n{#1}_r}}  
\ncomm{\Li}[1]{\em{{\n{#1}_{\infty}}}}
\ncomm{\Lm}[1]{\em{{\n{#1}_{\infty}}}}

\ncomm{\No}[1]{\em{\n{#1}_1}}
\ncomm{\Nt}[1]{\em{{\n{#1}}_2}}
\ncomm{\Np}[1]{\em{\n{#1}_p}}
\ncomm{\Nq}[1]{\em{\n{#1}_q}}
\ncomm{\Nr}[1]{\em{\n{#1}_r}}  
\ncomm{\Ni}[1]{\em{{\n{#1}_{\infty}}}}
\ncomm{\Nm}[1]{\em{{\n{#1}_{\infty}}}}

\ncomm{\Ltt}[1]{\em{{\n{#1}}_2^2}}

\def\ol{\overline}

\def\lam{\lambda}

\def\SM{\!\setminus\!}
\ncomm{\BS}[1]{_{{#1}\backslash}}
\ncomm{\FS}[1]{_{\slash #1}}
\ncomm{\doframe}[2]{
	\doit{frame}{
		\frametitle{#1}
		#2
	}
}

\ncomm{\doblock}[2]{
	\doit{block}{{#1}{#2}}
}

{\catcode`\|=\active
  \gdef\Set#1{\left\{\:{\mathcode`\|"8000\let|\SetVert #1}\:\right\}}}
{\catcode`\|=\active
  \gdef\Pres#1{\left\langle\:{\mathcode`\|"8000\let|\SetVert #1}\:\right\rangle}}
\def\SetVert{\egroup\;\middle|\;\bgroup}

\ncomm{\hide}[1]{}
\ncomm{\showit}[1]{#1}

\NCOM{\doMINI}[1]{
	\doit{minipage}{{\textwidth}#1}
}

\NCOM{\addTXT}[1]{
		\item[]
		\hskip-\leftmargin 
		\doMINI{#1}
		\bigskip
}

\NCOM{\doBIB}[2]{\doit{thebibliography}{{#1} #2}}

\NCOM{\doBIBtxt}[3]{
	\doit{thebibliography}{{#1} 
		\addTXT{#2}
		#3
	}
}

\NCOM{\MOD}[1]{(\tx{mod }#1)}

\NCOM{\useBIB}[2]{\bibliographystyle{#1}\bibliography{#2}}
\NCOM{\addBIB}[1]{\bibliography{#1}}

\NCOM{\EQ}[1]{Equation \re{#1}}